\documentclass[12pt,reqno]{article}
mm \textheight=8.9in

\usepackage{amssymb}
\usepackage{amsmath}

\usepackage{amsthm}
\usepackage{graphics}
\usepackage{color}
\newcommand{\tw}[3]{{$#1$}${\,\scriptscriptstyle {#2}}\atop\raise9pt\hbox{$\scriptstyle\tp$} ${$#3$}}
\newcommand{\st}[1]{\mbox{${\,\scriptscriptstyle {#1}}\atop\raise5.5pt\hbox{$*$}$}}
\newcommand{\btr}{\raise1.2pt\hbox{$\scriptstyle\blacktriangleright$}\hspace{2pt}}

\newcommand{\id}{\mathrm{id}}

\newcommand{\Lc}{\mathcal{L}}
\newcommand{\A}{\mathcal{A}}

\newcommand{\z}{\mathfrak{z}}

\newcommand{\Ha}{\mathcal{H}}
\newcommand{\Ru}{\mathcal{R}}

\newcommand{\Q}{\mathbb{Q}}

\newcommand{\C}{\mathbb{C}}
\newcommand{\Z}{\mathbb{Z}}

\newcommand{\N}{\mathbb{N}}

\newcommand{\tp}{\otimes}

\newcommand{\zt}{\zeta}

\newcommand{\la}{\lambda}

\newcommand{\End}{\mathrm{End}}

\newcommand{\Span}{\mathrm{Span}}

\newcommand{\btl}{\mbox{\raise1.1pt\hbox{$\scriptstyle\blacktriangleleft$}}}

\newcommand{\n}{\mathfrak{n}}

\newcommand{\nn}{\nonumber}

\newcommand{\K}{\mathcal{K}}

\newcommand{\si}{\sigma}
\newcommand{\al}{\alpha}

\newcommand{\bt}{\beta}

\newcommand{\be}{\begin{eqnarray}}
\newcommand{\ee}{\end{eqnarray}}

\newtheorem{thm}{Theorem}[section]
\newtheorem{propn}[thm]{Proposition}
\newtheorem{lemma}[thm]{Lemma}
\newtheorem{corollary}[thm]{Corollary}

\theoremstyle{definition}
\newtheorem{remark}[thm]{Remark}
\newtheorem{remarks}[thm]{Remarks}

\begin{document}
\title{Baxterization of solutions to reflection equation with Hecke R-matrix\footnote{
This research is supported by EPSRC grant C511166,
partially supported by the Minerva Foundation of Germany, and by
the RFBR grant no. 03-01-00593.} }
\author{P. P. Kulish$^\dag$ and A. I. Mudrov$^{\dag,\ddag,\natural}$}
\date{}
\maketitle \hspace{20pt}
\vbox{\noindent \small
$^\dag$St.-Petersburg Department of Steklov Mathematical
Institute, Fontanka 27,
\\
\vbox{191023 St.-Petersburg, Russia.
}
\\
$^{\ddag}$Department of Mathematics, University of York, YO10 5DD, UK,\\
$^{\natural}$Emmy Noether Institute for Mathematics, 52900 Ramat Gan, Israel,\\
\vbox{ e-mail: \texttt{kulish@pdmi.ras.ru, mudrov@pdmi.ras.ru}}}
\begin{abstract}
Let $R$ be a  Hecke solution to the Yang-Baxter equation
and $K$ be a reflection equation matrix with coefficients in an associative algebra $\A$.
Let $R(x)$ be the baxterization of $R$ and suppose that $K$ satisfies a
polynomial equation with coefficients in the center of $\A$.
We construct solutions to the
reflection equation with spectral parameter relative to $R(x)$,
in the form of polynomials in $K$.
\end{abstract} {\small \underline{Key words}: Reflection equation, baxterization.
\section{Introduction}
In this paper, we construct baxterization of reflection equation matrices associated
with Hecke solutions to the Yang-Baxter equation. It is assumed that the reflection
equation matrix satisfies a polynomial equation of degree $n+1$ with coefficients in some
commutative algebra.
The simplest non-trivial case $n=1$ was considered in \cite{FGX}.
A relevant problem was studied also in \cite{LM} using extensions of the affine Hecke algebra.
Under certain conditions, those extensions admit homomorphisms to the affine Hecke algebra,
yielding  solutions to the baxterization problem.

In the present paper we pursue a different approach, as compared to \cite{LM},
working directly within the affine Hecke algebra.
We give an alternative construction of reflection equation (RE) matrices with spectral parameter, using different basis.
Our choice of basis turns out to be very natural, leading to simple
formulas for the baxterization.

Let us describe the problem in more detail.
Fix  a field $k$ and let $V$ be a finite dimensional vector space over $k$.
An operator $R\in \End(V^{\tp 2})$ is called R-matrix if it satisfies the Yang-Baxter
equation
\be
R_{12}R_{23}R_{12}=R_{23}R_{12}R_{23}
\ee
in $\End(V^{\tp 3})$.
The lower indices specify the way of embedding $\End(V^{\tp 2}) \to\End(V^{\tp 3})$ in the usual way.
In this paper, we will deal with R-matrices satisfying the Hecke condition
\be
R^2=\omega R+1,
\label{Hecke}
\ee
where $\omega:=q-q^{-1}$ for some invertible $q\in k$ and $q\not=\pm 1$.

It is known that the rational function $R(x):=R-x^{-1}R^{-1} \in \End(V^{\tp 2})(x)$ satisfies the Yang-Baxter equation
with spectral parameter,
\be
R_{23}(x)R_{12}(xy)R_{23}(y)=R_{12}(y)R_{23}(x y)R_{12}(x),
\ee
in $\End(V^{\tp 2})(x,y)$.

Let $\A$ be an associative unital algebra over $k$.
A matrix $K\in \End(V)\tp \A$ is said to satisfy (constant) reflection equation (with coefficients in $\A$)
if
\be
RK_2RK_2=K_2RK_2R.
\label{re0}
\ee
We call $K$ an RE matrix with coefficients in $\A$.
Equation (\ref{re0}) is supported in $\End(V^{\tp 2})\tp \A$, and $K_2$ stands for $1\tp K$.
In particular, one can consider the RE algebra, which is the quotient of the free algebra $k\langle K^i_j\rangle$
modulo the relations (\ref{re0}). Then an RE matrix
with coefficients in an associative algebra $\A$ specifies a homomorphism
from the RE algebra to $\A$.

We will consider the following spectral version of the equation (\ref{re0}):
\be
R(\frac{x}{y})K_2(x)R(xy)K_2(y)=
K_2(y)R(xy)K_2(x)R(\frac{x}{y}).
\label{REspec}
\ee
This equation is supported in $\End(V^{\tp 2})\tp \A(x,y)$, with $K(x)$ being a rational function taking values in
$\End(V)\tp \A$. Equation (\ref{REspec}) was introduced in \cite{Cher1} to describe
the motion of a particle on half-line. It was applied to
integrable spin chains with non-periodic boundary conditions in \cite{S}. More references
to applications can be found in \cite{Do}.

We will consider the following problem: given a solution $K$ to (\ref{re0}), construct
a solution $K(x)$ to (\ref{REspec}) as a function of $K$. This procedure is called
baxterization of $K$. We solve this problem for a Hecke
R-matrix $R$ and $K$ satisfying a polynomial equation with coefficients in the center of $\A$.
The matrix $K(x)$ constructed in the present paper is
a rational function in $x$ and polynomial in $K$. Since the reflection
equation is invariant under multiplication by a scalar
function, $K(x)$ can be made polynomial
in $x$.

It is natural to consider the problem of baxterization in the abstract way, that is,
within the affine Hecke algebra $\Ha$ of rank one. The latter is the quotient of  $\z\langle R,K\rangle $ by the ideal of
relations (\ref{Hecke}) and (\ref{re0}).  In the algebra $\Ha$ we impose the relation
$p(K)=0$, where $p$ is a polynomial with coefficients in $\z$. The corresponding
quotient denoted further by $\Ha_p$ is known as the cyclotomic Hecke algebra, see e.g. \cite{Ar}.
We construct solutions of the equation (\ref{REspec}) in $\Ha_p$
and prove that there are no other solutions if $\z$ is a a field
and $p(0)\not=0$. The latter condition means invertibility
of $K$ in $\Ha_p$.

The paper is arranged as follows. In Section \ref{sRKa} we replace
the problem by a corresponding problem in the the affine Hecke
algebra $\Ha$. We reduce the equation (\ref{REspec}) to a certain linear subspace in
$\Ha$. In Section \ref{sSRE} we construct the baxterization.
In Section \ref{s2T} we show that the list of solutions from
Section \ref{sSRE} is exhaustive, under certain assumptions.
In Section \ref{secA} we consider some applications to the quantum group
$U_q\bigl(gl(n+1)\bigr)$.

\vspace{20pt}
\noindent
{\large\bf Acknowledgment.} We thank  I. Cherednik and V. Tarasov for
interesting discussions. P.K. is grateful to
the Laboratoire Physique Theorique et Modelization of
University Cergy-Pontoise for hospitality and to CNRS for financial
support.
\section{Baxterization in  the affine Hecke-algebra}
\label{sRKa}
\subsection{Affine Hecke-algebra of rank one}
We start with recalling the definition and basic properties of
the affine Hecke algebra of rank one.

We fix a ground field $k$ of zero characteristic and a commutative unital $k$-algebra $\z$.
The affine Hecke algebra $\Ha$ is the quotient of the free algebra
$\z\langle R,K\rangle$ by the ideal of relations
\be
R^2&=&\omega R+1\quad \mbox{(the Hecke condition)},
\label{HC}\\
RKRK&=&KRKR\quad \mbox{(the reflection equation)}.
\label{re1}
\ee
Here $\omega:=q-q^{-1}$ for some invertible $q\in k$ and $q\not=\pm 1$.

Put $Y:=RKR$; the reflection equation (\ref{re1}) implies that $Y$  commutes with $K$.
The following proposition is called PBW theorem for $\Ha$.
\begin{propn}[\cite{Cher2}]
\label{PBW}
The monomials $K^iY^jR^m$, where $i$ and $j$ are non-negative integers
and $m=0,1$, form a basis in $\Ha$.
\end{propn}
Let us introduce the $\z$-submodule $\Lc\subset \Ha$
spanned by the elements
$$e_{i;j}:=K^iRK^j-K^jRK^i,\quad i,j=0,1,\ldots.$$
Clearly $e_{i;j}=-e_{j;i}$.
\begin{lemma}
The elements $e_{i;j}$, where $i>j\geq 0$, form a basis in $\Lc$.
\end{lemma}
\begin{proof}
All we need to check is that the elements $e_{i;j}$ are independent
for different pairs of $i,j$ such that $i>j\geq 0$.
Using double induction, one can check the equality
\be
K^iRK^j-K^jRK^i&=&(K^iY^j-K^jY^i) R-\omega Y\frac{K^iY^{j}-K^jY^{i}}{K-Y}.
\label{KK-KY}
\ee
Now the lemma follows from Proposition \ref{PBW}.
\end{proof}

The submodule $\Lc$ plays an important role in our exposition
due to the following fact.
\begin{lemma}
\label{+-}
In the algebra $\Ha$, the following relation
holds for all $m\in\Z$ and all $k\in \N$:
\be
RK^mRK^k-K^kRK^mR=\omega\sum_{i=1}^k( K^{m+k-i}RK^i-K^iRK^{m+k-i}).
\label{+}
\ee
Let $\bar\Ha$ denote the quotient of the algebra $\Ha[\bar K]$,
by the ideal of relations $\bar K K=K\bar K=1$. Then for all $m\in\Z$ and all $k\in \N$
\be
RK^mRK^{-k}-K^{-k}RK^mR=\omega\sum_{i=1}^k(K^{i-k}RK^{m-i}-K^{m-i}RK^{i-k})
\label{-}
\ee
in $\bar\Ha$.
\end{lemma}
\begin{proof}
First let us check the formula (\ref{+}) for $k=1$. The general case will follow by induction on $k$
using the standard formula for the commutator derivation in associative algebras.

Wherever it applies, we assume $K$ to be invertible.
For any $m\in \Z$ the identity (\ref{re1}) implies the equality $K^mRKR=RKRK^m$ and therefore
the equality
\be
RK^mRKRKR=RKRKRK^mR.
\label{aux_eq1}\ee
We multiply both sides of the equation (\ref{re1}) by $RK^{m-1}$ from the left and by $R$ from the right.
Then we multiply (\ref{re1}) by $K^{m-1}R$ from the right and by $R$ from the left.
That gives the following two equalities:
$$
\left\{
\begin{array}{lcc}
RK^{m}RKR^2=RK^{m-1}RKRKR, \nn\\
RKRKRK^{m-1}R=R^2KRK^{m}R. \nn
\end{array}
\right.
$$
Taking into account (\ref{aux_eq1}) and the Hecke condition (\ref{HC}), from this we deduce
\be
RK^{m}RK(\omega R+1)&=&(\omega R+1)KRK^{m}R
\quad\mbox{or}\nn
\\
RK^{m}RK-KRK^{m}R&=&\omega (RKRK^{m}R-RK^{m}RK R)\nn\\
&=&\omega (K^{m}RKR^2-R^2K RK^{m})\nn\\
&=&\omega (K^{m}RK-K RK^{m})\nn
\ee
for any $m\in \Z$. Here we have subsequently used the reflection equation and
the Hecke condition.
The formula (\ref{+}) is thus proven for $k=1$. The general case $n>0$ follows by induction through
the commutator derivation
\be
[RK^{m}R,K^{k}]&=&K[RK^{m}R,K^{k-1}]+[RK^{m}R,K]K^{k-1}.\nn
\ee
The formula (\ref{-}) is obtained from (\ref{+}) via the left and right multiplication by $K^{-k}$.
\end{proof}

\subsection{Reflection equation with spectral parameter}
Put $R(x):=R-x^{-1}R^{-1}\in \Ha(x)=\Ha\tp k(x)$.
In the algebra  $\Ha(x,y)$, consider the equation
\be
\label{re}
R(\frac{x}{y})K(x)R(xy)K(y)=
K(y)R(xy)K(x)R(\frac{x}{y})
\ee
with respect to $K(x)\in \Ha(x)$. The problem is to construct solutions to
this equation as functions of $K$.
For any element $f(x)\in \z(x)$
the scalar matrix $f(x)\id_V$ clearly
solves (\ref{re}). Such solutions are called trivial.
In general, if $K(x)$ is a solution, then $f(x)K(x)$ is a solution as well.

The following is the key observation in the theory of (\ref{re}).
\begin{propn}
The equation {(\ref{re})} reduces to an equation in $\Lc(x,y)=\Lc\tp k(x,y)$.
\label{reduction}
\end{propn}
\begin{proof}
The Hecke condition (\ref{HC}) implies $R(x)=R(1-x^{-1})+\omega x^{-1}$.
Substitute this into {(\ref{re})} and rewrite the latter as
\be
&(x-\frac{1}{y}-y+\frac{1}{x})[R K(x)RK(y)-K(y)R K(x)R]+
\nn\\
&+\omega(y-\frac{1}{x})[K(x) R K(y)-K(y)R K(x)]
+\omega(\frac{1}{y}-\frac{1}{x})[R K(x) K(y)- K(y)  K(x) R]=0.
\label{re'}
\ee
The second and third  difference terms obviously lie in $\Lc(x,y)$.
To prove the proposition we must show the same for the first difference.
That follows from Lemma \ref{+-}.
\end{proof}

\subsection{Cyclotomic Hecke algebra}
We will solve the equation (\ref{re}) imposing a polynomial condition on the
generator $K$. This subsection is devoted
to basic properties of the corresponding quotient of
the algebra $\Ha$.

Suppose that $n\geq 1$ and choose a sequence $(A_i)_{i=0}^n$ of elements from $\z$.
Define a polynomial $p$ in one variable setting
\be
p(K)=-K^{n+1}+A_0K^n+A_1K^{n-1}+\ldots+A_n.
\label{polynomial}
\ee
We will need the following elementary fact, which is true for any polynomial
with invertible highest coefficient.
\begin{lemma}
The quotient algebra $\z[K]/(p)$ is a free $\z$-module generated by
 $\{K^i\}_{i=0}^\n$.
\label{free0}
\end{lemma}

Let $\K$ denote the subalgebra in $\Ha$ generated by $K$. It follows from Proposition \ref{PBW}
that $\K$ is isomorphic to $\z[K]$.
Consider  the quotient $\Ha_p$ of the algebra $\Ha$ by the ideal generated by $p(K)$.
This quotient is called {\em cyclotomic Hecke algebra}, \cite{Ar}, under the assumption
that $p$ is decomposable in $\z[K]$ by linear divisors.
By $\K_p\subset \Ha_p$ we denote the corresponding quotient of the subalgebra $\K$.
Proposition \ref{PBW} implies that $\K_p\simeq \K/(p)$; by Lemma \ref{free0}
$\K_p$ is free over $\z$.

Define polynomials $p_\pm$ in two variables setting $p_\pm(K,Y)=\pm q^{\pm 1}p(Y)-\omega Y\frac{p(Y)-p(K)}{Y-K}$.
Let $J_\pm$ denote the ideals in $\z[K,Y]$ generated by
$p_\pm(K,Y)$ and by $p(K)$.
\begin{lemma}
As a $\z$-module, the algebra $\Ha_p$ is isomorphic to $\z[K,Y]/J_+\oplus \z[K,Y]/J_-$.
\label{free1}
\end{lemma}
\begin{proof}
The $k$-algebra generated by $R$ is two dimensional semisimple.
Let $P_\pm$ denote the idempotents $P_\pm:=\frac{\pm q^{\pm 1}-R}{q+q^{-1}}$.
Then $\Ha$ splits into the direct sum of $\z$-modules,
$\z[K,Y]P_+\oplus \z[K,Y]P_-$.
Now the statement follows from the formula
$R p(K)=p(Y)R-\omega Y\frac{p(Y)-p(K)}{Y-K}$.
\end{proof}

\begin{propn}
\label{free}
The algebra $\Ha_p$ is a free $\z$-module generated by $K^iY^jR^m$, where $i,j=0,\ldots,n$, $m=0,1$.
\end{propn}
\begin{proof}
In view of Lemma \ref{free1}, it suffices to proof that both algebras $\z[K,Y]/J_\pm$
are freely generated over $\z$ by $K^iY^j$, where $i,j=0,\ldots,n$. Now the proposition
follows from Lemma \ref{free0}. Indeed, consider $\z[K,Y]/J_\pm$ as
the quotient of $\K_p[Y]$ by the image of the ideal $(p_\pm)$ under the projection
$\z[K,Y]\to \K_p[Y]$. The highest coefficients
of $p_\pm$ as polynomials in $Y$ are $\pm q^{\mp 1}$ and hence invertible. Therefore $\z[K,Y]/J_\pm$
are freely generated over $\K_p$
by $\{Y^j\}_{j=0}^n$. On the other hand,  $\K_p$ is freely generated
over $\z$ by $\{K^j\}_{j=0}^n$. This proves the statement.
\end{proof}
Let $\Lc_p$ denote the image of the submodule $\Lc\in \Ha$ under the projection $\Ha\to \Ha_p$.
\begin{corollary}
As a $\z$-module, $\Lc_p$ is freely generated by $\{e_{j;i}\}_{0\leq i<j\leq n}$.
\label{L_free}
\end{corollary}
\begin{proof}
The elements $\{e_{j;i}\}_{0\leq i<j\leq n}$ clearly generate $\Lc_p$ over $\z$.
Their independence follows from Proposition \ref{free} and formula (\ref{KK-KY}).
\end{proof}

\section{Solving RE with spectral parameter}
\label{sSRE}
In this section we guess a solution to (\ref{re})
and prove that it satisfies (\ref{re}) indeed. In the subsequent section
we show that essentially there are no other solutions if $\z$ is an integral domain (has no zero divisors) and
$p(0)\not=0$.
These conditions are
fulfilled e.g. for $\z$ being the center of the reflection equation algebra
associated with $U_q\bigl(gl(n+1)\bigr)$, cf. Section \ref{secA}.

Let us introduce the dual basis $\{\al^{(n)}_{j;i}\}_{0\leq i<j\leq n}\subset \Lc_p^*$
to $\{e_{j;i}\}_{0\leq i<j\leq n}$.
We view the left-hand side of (\ref{re'}) as a quadratic map
$\K(x)\to \Lc(x,y)$. This map induces
a map
$\K_p(x)\to \Lc_p(x,y)$, which will be further denoted by $Q_p$.
In view of Proposition \ref{reduction}, Corollary  \ref{L_free} allows to
reduce the equation (\ref{re}) to the system of functional equations
\be
(\al^{(n)}_{j;i}\circ Q_p)\bigl(K(x)\bigr)=0, \quad 0\leq i<j \leq n
\label{system}.
\ee

For $n>1$ consider the element $\tilde K:=K^n-A_0K^{n-1}-\ldots-A_{n-1}$; it satisfies
the equality $K\tilde K=A_n$. Obviously the elements $\{K^i\}_{i=0}^{n-1}\cup \{\tilde K\}$
form a basis in $\K_p$. We will call a solution to (\ref{re}) {\em principal} if either $n=1$ or $n>1$ and
the coefficient before $K^{n-1}$ in the expansion over this basis is non-zero. Otherwise, the solution will  be called
 {\em small}.
\subsection{Example: quadratic relation on $K$}
Let us start from the case $n=1$.
Impose on the matrix $K$ the quadratic equality
$K^2=A_0K+A_1$.
Then the subspace $\Lc_p$ is one dimensional and spanned by the element $e_{1;0}$.
We get the only equation on the coefficients $a_0(x)$ and $a_1(x)$ in the expansion $K(x)=a_0(x)K+a_1(x)$:
\be
(x-\frac{1}{y}-y+\frac{1}{x})a_1'a_0''
+(y-\frac{1}{x})(a_1' a_0''-a_0'a_1'')
+(\frac{1}{y}-\frac{1}{x})(a_0'a_0''A_0+a_0'a_1''+a_1'a_0'')=0.
\nn
\ee
Here $a'_i:=a_i(x)$ and $a''_i:=a_i(y)$.
Solving this equation we recover the following result.
\begin{propn}[\cite{FGX}]
\label{FGX}
Let $K$ be the generator of the algebra $\K_p$, with $p(K)=-K^2+A_0K+A_1$, $A_i\in \z$.
Then the element
\be
K(x):=K+\frac{\zt-x^{-1}(\xi x- A_0)}{(x-x^{-1})}\in \K_p(x),
\label{quadratic}
\ee
where $\zt$ and $\xi$ are arbitrary elements from $\z$,
solves the equation (\ref{re}).
If $\z$ is a field, then any non-trivial solution is given by
(\ref{quadratic}) up to a factor from $\z(x)$.
\end{propn}
Remark that although the right-hand side of (\ref{quadratic}) depends only
on the difference $\zt-\xi$, it is convenient to retain both parameters,
to unify (\ref{quadratic}) with the general case $n\geq 2$.
\subsection{Example: cubic relation on $K$}
When $n=2$ the space $\Lc_p$ is spanned by the
elements $e_{2;1}$, $e_{2;0}$, and $e_{1;0}$.
An element $K(x)\in \Ha_p(x)$ is expanded as
$K(x)=a_0(x)K^2+a_1(x)K+a_2(x)$, where
$a_i(x)$ are some rational functions of  $x$ with values in $\z$.

Equation (\ref{re}) is equivalent to the system of functional equations
$$
\left\{
\begin{array}{lcc}
(x-\frac{1}{y}-y+\frac{1}{x})a_0'a_0''A_0
+(x-\frac{1}{y})a_0'a_1''
-(y-\frac{1}{x})a_0''a_1'&=&0,\nn
\\[4pt]
(x-\frac{1}{x})a_0''a_2'
-(y-\frac{1}{y})a_0'a_2''-
(\frac{1}{y}-\frac{1}{x})\Bigl(a_0'a_0''(A_0^2+A_1)+\bigl(a_0'a_1''+a_0''a_1'\bigr)A_0
+a_1'a_1''\Bigr)&=&0,
\\[4pt]
(x-y)A_2a_0'a_0''
+(x-\frac{1}{x})a_1''a_2'
-(y-\frac{1}{y})a_1'a_2''
+(\frac{1}{y}-\frac{1}{x})\Bigl( \bigl(a_0'a_1''+a_0''a_1'\bigr)A_1
+a_0'a_0''A_0A_1\Bigr)
&=&0.
\end{array}
\right.
$$
Here $a'_i:=a_i(x)$ and $a''_i:=a_i(y)$.

This system is easy to solve, and the result is given by the following proposition.
\begin{propn}
\label{thm_cubic}
Let $K$ be the generator of the algebra $\K_p$, with $p(K)=-K^3+(A_0K^2+A_1K+A_2)$, $A_i\in \z$.
For any pair $\xi, \zt\in \z$ such that $\zeta \xi=A_2$
the element
\be
K(x):=K^2+(\xi x-A_0)K+\frac{\zeta-x^{-1}(x^2\xi^2-A_0\xi x-A_1)}{(x-x^{-1})}\in \K_p(x),
\label{cubic}
\ee
solves (\ref{re}).
If $\z$ is a field, then any non-trivial solution is given by
(\ref{cubic}) up to a factor from $\z(x)$.
\end{propn}
Note that the solution (\ref{cubic}) is principal if $\xi\not= 0$
and small if $\xi=0$. The latter case is possible only if $A_2=0$.
\subsection{Finding small solution}
\label{subsecSmall}
In this subsection we assume that
$K$ is invertible in $\K_p\subset \Ha_p$.
With $\z$ being an integral domain and invertible $K$,
Corollary \ref{small} below states that a small solution to (\ref{re}) may exist only if $n= 3$.
Then it can be presented as a linear combination $a_+K+a_0+a_-K^{-1}$, where
$a_\pm$ and $a_0$ are some rational functions with values in $\z$.
In this section we will find the small solution explicitly.

We assume that element $-A_3=-p(0)$ has a square root $\sqrt{-A_3}$
in $\z$.
\begin{propn}
Suppose that $K$ is invertible in $\K_p$, with
$p(K)=-K^{n+1}+\sum_{i=0}^n A_i K^{n-i}$.
A non-trivial  small solution  to (\ref{re})
exists only for $n=3$ and has the form
\be
\label{small_sol}
K(x)=x\sqrt{-A_3} K+\frac{A_0\sqrt{-A_3}+A_2x}{x-x^{-1}}+A_3K^{-1}
\in \K_p(x).
\ee
This small solution is unique up to a factor from $\z(x)$,
provided $\z$ is an integral domain.
\end{propn}

\begin{proof}
We will work in the basis $\{K^i\}_{i=-2}^1$. Then the submodule $\Lc_p$ is generated over $\z$ by
the elements $e_{i;j}:=K^iRK^j-K^jRK^i$, where $i,j=-2,\ldots, 1$, $i>j$.
Expanding  (\ref{re'}) over this basis, we get for (\ref{small_sol}) the
 following system of  functional equations:
\be
\left\{
\begin{array}{lcc}
-(x-\frac{1}{y}-y+\frac{1}{x})a_-'a_+''+(y-\frac{1}{x})(a_+'a_-''-a_+''a_-')&=&0,\\
(x-\frac{1}{y}-y+\frac{1}{x})a_-'a_-''+(\frac{1}{y}-\frac{1}{x})(a_-'a_-''+A_3a_+''a_+')&=&0,\\
(x-\frac{1}{y}-y+\frac{1}{x})a_0'a_-''+(y-\frac{1}{x})(a_0'a_-''-a_0''a_-')
+(\frac{1}{y}-\frac{1}{x})(a_-'a_0''+a_-''a_0'+A_2a_+''a_+')&=&0,\\
(x-\frac{1}{y}-y+\frac{1}{x})a_0'a_+''
-(y-\frac{1}{x})(a_+'a_0''-a_+''a_0')+(\frac{1}{y}-\frac{1}{x})(a_+'a_0''+a_+''a_0'+A_0a_+''a_+')
&=&0.
\end{array}
\right.
\nn
\ee
If $a_+= 0$, then $a_-= 0$ as well, and such a solution is trivial.
Suppose that $a_+\not = 0$. Setting $a_+(x)=\sqrt{-A_3}\>x\not=0$, we continue
\be
\left\{
\begin{array}{lcc}
-a_-'+a_-''&=&0,\\
a_-'a_-''-A_3^2&=&0,\\
(x-\frac{1}{x})a_0'a_-''-(y-\frac{1}{y})a_0''a_-'
-(x-y)A_2A_3&=&0,\\
(x-\frac{1}{x})ya_0'
+(y-\frac{1}{y})xa_0''+xy(\frac{1}{y}-\frac{1}{x})A_0\sqrt{-A_3}
&=&0,
\end{array}
\right.
\nn
\Leftrightarrow
\left\{
\begin{array}{lcc}
a_-&=&A_3,\\
\zt &=& A_2,\\
(x-\frac{1}{x})a_0
&=&
(A_0\sqrt{-A_3}+x\zt).
\end{array}
\right.
\nn
\ee
This shows that (\ref{small_sol}) is a solution to (\ref{re}).
This also proves its uniqueness if $\z$ is an integral domain.
\end{proof}
Thus there are two different small solutions corresponding to
the two different square roots of $-A_3$ in $\z$.

In the basis
$\{K^i\}_{i=0}^3$ the element (\ref{small_sol}) can be presented as
\be
K(x)&=&
K^3-A_0K^2
-(A_1-x\sqrt{-A_3})K+\frac{x^{-1} A_2+ A_0\sqrt{-A_3}}{x-x^{-1}}.\nn
\ee
In the limit $A_3\to 0$ this coincides
with the limiting case $A_3=0$, $\xi=0$ of the solution (\ref{main}).

\subsection{Preparatory technical material}
In this subsection we develop a certain machinery to finding principal solutions to (\ref{re}).

For $i=0,\ldots, n+1$ we define the polynomials $\tilde \phi_i\in \z[u]$ to be the
regular parts of the Laurent polynomials $p(u)u^{i-(n+1)}\in \z[u,u^{-1}]$, where $p$ is given by (\ref{polynomial}).
In particular, $\tilde \phi_{n+1}(u)=p(u)$.
Introduce the functions $\phi_i(\xi,x)$ setting $\phi_i(\xi,x):=-\tilde\phi_i(\xi x)$.
We will need the following recurrent formulas
\be
\phi_{0}(\xi, x):=1
,\quad
\phi_{i+1}(\xi, x):=\phi_{i}(\xi, x)\xi x-A_{i},
\label{phi}
\ee
where $A_i$ are the coefficients of the polynomial $p$.
The functions $\phi_{k}(\xi, x)$ will play an important role in our exposition.

Below we derive some identities for $\phi_{k}$ that we employ later on.
To simplify the  formulas  we use notation $a':=a(x)$ and $a'':=a(y)$
for an element  $a\in \z[x]$. Put also $\Delta_i:=\phi_i'-\phi_i''$ for $i\in [0,n+1]$.
It follows, in particular, that $\Delta_0=0$.
\begin{lemma}
Let $i$ and $m$ be non-negative integers such that $0\leq i\leq m\leq n$. Then
\be
\xi(x-y)\sum_{\al=i}^{m} \phi_\al'\phi_{m+i-\al}''
=
\phi_{m+1}'\phi_{i}''-\phi_{i}'\phi_{m+1}''
-\sum_{\al=i}^{m}A_\al \Delta_{m+i-\al}.
\label{eq_dif}
\ee
\label{dif}
\end{lemma}
\begin{proof}
We have for the right-hand side of (\ref{eq_dif}):
\be
\mbox{r.h.s. of } (\ref{eq_dif})&=&
\sum_{\al=i}^{m} \xi x\phi_\al'\phi_{m+i-\al}''-\sum_{\al=i}^{m} \phi_\al'\xi y\phi_{m+i-\al}''
=\nn\\
&=&\sum_{\al=i}^{m}(\phi_{\al+1}'+\phi_0'A_\al)\phi_{m+i-\al}''-\sum_{\al=i}^{m} \phi_{m+i-\al}'(\phi_{\al+1}''+\phi_0''A_\al)
=\nn\\
&=&\phi_{m+1}'\phi_{i}''-\phi_{i}'\phi_{m+1}''+
(\ldots)
-\sum_{\al=i}^{m}A_\al \Delta_{m+i-\al}
\nn
\ee
The term in the parentheses in the bottom line equals
$\sum_{\al=i}^{m-1}\phi_{\al+1}'\phi_{m+i-\al}''-\sum_{\al=i}^{m-1} \phi_{m+i-\al}'\phi_{\al+1}''=0$.
This proves the lemma.
\end{proof}

For each $m\in [n+1,2n]$ we introduce the elements $\{A^{(m)}_i\}_{i=0}^n\subset \z$ to be the
coefficients of the expansion
$$
K^m=\sum_{i=0}^n A^{(m)}_i K^{n-i}\in \K_p.
$$
By definition, $A^{(n+1)}_i=A_i$, i.e. the coefficients of the polynomial $p$.

The following rule of sums holds true for $\{A^{(m)}_i\}$.
\begin{lemma}
Let $i,k$ be non-negative integers such that  $0\leq i\leq n$ and $0\leq k<n$. Then
\be
A^{(n+k+1)}_i&=&\sum_{\al=0}^{k-1}A^{(n+k-\al)}_{i}A_\al+A_{k+i},\quad i=0,\ldots, n-k,
\label{aux_eq5}
\\
A^{(n+k+1)}_i&=&\sum_{\al=0}^{k-1}A^{(n+k-\al)}_{i}A_\al,\quad i=n-k+1,\ldots, n.
\label{aux_eq6}
\ee
\label{sums}
\end{lemma}
\begin{proof}
The lemma immediately follows from the presentation
\be
K^{n+k+1}&=&
K^{k}(A_0K^n+A_1K^{n-1}+\ldots+A_{n-1}K^1+A_n)\nn\\\
&=&(A_0K^{n+k}+A_1K^{n-1+k}+\ldots+A_{k-1}K^{1+n})+(A_k K^{n}+\ldots+A_{n-1}K^{1+k}+A_nK^k).
\nn
\ee
The first summand contributes to the r.h.s. of
(\ref{aux_eq5}) and (\ref{aux_eq6}), while the second one only to (\ref{aux_eq5}).
\end{proof}
\begin{lemma}
Let $i$ and $j$ be non-negative integers such that $0<i\leq n$ and  $0\leq j\leq n$. Then
\be
\xi(x-y)\sum_{\al=0}^{i-1} A^{(n+i-\al)}_j\sum_{\bt=0}^{\al}\phi_\bt'\phi_{\al-\bt}''
&=&\sum_{\al=\max\{j+i-n,0\}}^{i}\Delta_\al A_{j+i-\al}.
\label{aux_eq4}
\ee
\label{sums'}
\end{lemma}
\begin{proof}
Applying Lemma \ref{dif} we find
\be
\mbox{l.h.s. of } (\ref{aux_eq4}) &=&
\sum_{\al=0}^{i-1} A^{(n+i-\al)}_j(\Delta_{\al+1}-\sum_{\bt=0}^{\al}A_{\al-\bt}\Delta_{\bt})
\nn\\
&=&A_j\Delta_{i}+
\sum_{\al=0}^{i-1}\Delta_\al( A^{(n+i+1-\al)}_j-\sum_{\bt=\al}^{i-1} A^{(n+i-\bt)}_jA_{\bt-\al})
\nn\\
&=&A_j\Delta_{i}+
\sum_{\al=0}^{i-1}\Delta_\al( A^{(n+(i-\al)+1)}_j-\sum_{\bt=0}^{i-1} A^{(n+(i-\al)-\bt)}_jA_{\bt})
\nn\\
&=&A_j\Delta_{i}+
\sum_{\al=\max\{j+i-n,0\}}^{i-1}\Delta_\al A_{j+i-\al}=\sum_{\al=\max\{j+i-n,0\}}^{i}\Delta_\al A_{j+i-\al},
\nn
\ee
as required.
Passing to the bottom line we used Lemma \ref{sums}.
\end{proof}
\subsection{The main theorem}
In the present section we show that a principal solution does exist and
has the explicit form as stated by the following theorem.
\begin{thm}
Fix $n\geq 2$. Let $K$ be the generator of the algebra $\K_p$,
with $p$ given by (\ref{polynomial}).
Then the element
\be
K(x):=\sum_{i=0}^{n-1}\phi_i(\xi, x)K^{n-i}+
\frac{\zt-x^{-1}\phi_n(\xi, x)}{x-x^{-1}}
\in \K_p(x),
\label{main}
\ee
where $\phi_i$ are defined by (\ref{phi}),
is a principal solution to (\ref{re})
for any pair $\xi, \zt\in \z$ such that $\zeta \xi=p(0)$.
\label{thm_main}
\end{thm}
\begin{proof}
Fix $\zt\in \z$ to be a divisor of $p(0)$. Zero $\zt$ is also admissible in the
case $p(0)=0$.
Consider the  polynomial $\tilde p$ in one variable with coefficients
in $\z[v]$ defined by $\tilde p(u):=-u^{n+1}+\sum_{i=0}^{n-1}A_iu^{n-i}+ v\zt$.
First let us prove the statement replacing the ring $\z$
by $\z[v]$ and the polynomial $p$ by $\tilde p$.
The parameter $\xi$ in $\phi_i(\xi,x)$ is replaced by the indeterminate $v$.
The element (\ref{main}) will have  polynomial dependence on $v$,
and we can consider (\ref{system}) as a system of polynomial equations with respect to $v$.
Therefore, to prove the statement for $\z[v]$, it is sufficient to
check it setting $v$ a non-zero element from $k$ (the field $k$ is
assumed to be infinite).

Equation (\ref{re}) is equivalent to the system (\ref{system}).
Fix a pair of non-negative integers $i,j$ satisfying the condition $0\leq i<j< n$.
The contribution of the first difference term  in (\ref{re'}) to the coefficient before $e_{n-i;n-j}$
equals
\be
(1-\frac{1}{xy})(x-y)\bigl(
\sum_{\al=0}^{j-1} A^{(n+j-\al)}_i\sum_{\bt=0}^{\al}\phi_\bt'\phi_{\al-\bt}''
-
\sum_{\al=0}^{i-1} A^{(n+i-\al)}_j\sum_{\bt=0}^{\al}\phi_\bt'\phi_{\al-\bt}''
+\sum_{\al=i}^{j}\phi_\al'\phi_{i+j-\al}''
-\phi_j'\phi_i''\bigr).
\nn
\ee
Using Lemmas \ref{dif} and \ref{sums'} we transform this expression to
\be
(1-\frac{1}{xy})\frac{1}{v}(-\Delta_i A_{j}+
\phi_{j+1}'\phi_{i}''-\phi_{j+1}''\phi_{i}'
)
-(x-\frac{1}{y}-y+\frac{1}{x})\phi_j'\phi_i''
.
\nn
\ee
The functional $\al^{(n)}_{n-i;n-j}$ vanishes on the last term in (\ref{re'}),
since $i,j$ are assumed to be less than $n$. The contribution of the second term
to the $e_{n-i;n-j}$-coefficient equals $(y-\frac{1}{x}) (\phi_i'\phi_j''-\phi_j'\phi_i'')$.
Thus we have
$$
(\al^{(n)}_{n-i;n-j}\circ Q_p)\bigl(K(x)\bigr)=(1-\frac{1}{xy})\frac{1}{v}(-\Delta_i A_{j}+
\phi_{j+1}'\phi_{i}''-\phi_{j+1}''\phi_{i}')+
(y-\frac{1}{x}) \phi_i'\phi_j''
-(x-\frac{1}{y}) \phi_j'\phi_i''.
$$
It is easy to see that this expression vanishes, taking into account $\Delta_i=\phi_i'-\phi_i''$ and
the definition (\ref{phi}) of functions $\phi_i$.

Now fix $i\in[0,n-1]$ and check that the functional $\al^{(n)}_{0;n-i}:=-\al^{(n)}_{n-i;0}$ vanishes on $Q_p\bigl(K(x)\bigr)$.
Denote by $b\in \z(x)$ the free term of (\ref{main}),
\be
b(x):=\frac{\zt-x^{-1}\phi_n(v, x)}{x-x^{-1}}.
\label{b}
\ee
A straightforward computation gives $(\al^{(n)}_{0;n-i}\circ Q_p)\bigl(K(x)\bigr)$
 to be equal to
\be
&(1-\frac{1}{xy})(x-y)(\sum_{\al=0}^{i-1} A^{(n+i-\al)}_n\sum_{\bt=0}^{\al}\phi_\bt'\phi_{\al-\bt}'')
+(x-\frac{1}{x})b'\phi_i''-(y-\frac{1}{y}) \phi_i'b''
\nn\\
&-(\frac{1}{y}-\frac{1}{x})(\phi_n'\phi_i'' +\phi_n''\phi_i' )+
(\frac{1}{xy})(x-y)(\sum_{\al=i}^{n} \phi_\al'\phi_{n+i-\al}''
+\sum_{j=0}^{n-1}A^{2n-j}_i\sum_{\al=0}^j \phi_\al'\phi_{j-\al}'').
\nn
\ee
Applying Lemmas \ref{dif} and \ref{sums'} we transform this expression to
\be
&(1-\frac{1}{xy})\frac{1}{v}(v\zt \Delta_i)
+(x-\frac{1}{x})b'\phi_i''-(y-\frac{1}{y}) \phi_i'b''
-(\frac{1}{y}-\frac{1}{x})(\phi_n'\phi_i'' +\phi_n''\phi_i' )+
\nn\\
&+
(\frac{1}{xy})\frac{1}{v}(\phi_{n+1}'\phi_{i}''-\phi_{i}'\phi_{n+1}'').
\nn
\ee
It is immediate to see that it vanishes once $b$ equals
(\ref{b}). Here we have used $\tilde p(0)=v\zt$.

Thus we have proven the theorem for the ring $\z[v]$ and the polynomial $\tilde p$.
Now if $\xi\in \z$ is such that $\xi\zt=p(0)$, the evaluation homomorphism
$\z[v]\to \z$, $v\mapsto \xi$, sends $\tilde p$ to $p$ and hence
extends to a homomorphism $\Ha_{\tilde p}\to \Ha_{p}$.
This proves the theorem for $\z$ and $p$.
\end{proof}
\begin{remarks}
\item
\begin{enumerate}
\item
Obviously, the set of pairs $\xi,\zt$ satisfying the condition
$\xi\zt=A_n$ is not empty. One can take, say for $\xi$,
an arbitrary non-zero element from $k$.
\item
If $\z$ is an integral domain, the coefficient $b$ is determined by (\ref{b}) uniquely.
Indeed, $b$ is specified by the equation $(\al^{(n)}_{0;n-i}\circ Q_p)\bigl(K(x)\bigr)=0$, $i=0$,
up to a constant summand from $ \z$. This constant is fixed to be $\zt$ such that $\zt\xi=A_n$,
by the equations corresponding to higher $i$. Such $\zt$ is unique if $\z$ has no zero divisors.
\item
Neither (\ref{small_sol}) nor (\ref{main}) depend on the parameter $q$
entering the definition of $\Ha$.
\end{enumerate}
\end{remarks}

\subsection{Properties of the principal solution}
Note that the formula (\ref{main}) gives a solution to (\ref{re}) for $n=1$ (and for $n=0$) as well.
The difference from (\ref{quadratic}) is that the parameters $\xi$ and $\zt$
in (\ref{main}) are bound by the condition $\xi\zt=p(0)$, while being
independent in (\ref{quadratic}).
Let $K_p(\xi,\zt,x)$ denote the principal solution (\ref{main}) for $n\geq 0$
If $\z$ is an integral domain, $\zt$ is determined by $\xi$ uniquely. We will hide
the dependance on $\zt$ if the latter is clear from the context.
Next we investigate links among $K_p(\xi,\zt,x)$ corresponding to different $p$.
\begin{propn}
\label{pp'}
Suppose that the polynomial $p$
 is divisible by $\tilde p$,
that is $p=\tilde{\tilde p} \tilde p$ for some polynomial $\tilde{\tilde p}$.
Let $\xi,\tilde \zt\in \z$ be such that $\xi\tilde \zt =\tilde p(0)$ and
put $\zt=\tilde{\tilde p}(0)\tilde \zt$.
Then the natural projection homomorphism $\theta\colon \Ha_{p}(x)\to \Ha_{\tilde p}(x)$
brings $K_{p}(\xi,\zt, x)$ to $\tilde{\tilde p}(\xi x)K_{\tilde p}(\xi,\tilde\zt ,x)$.
\end{propn}
\begin{proof}
One can directly check the statement for the case $\deg \tilde{\tilde p}=1$.
The obvious induction on $\deg \tilde{\tilde p}=\deg p-\deg \tilde p$ proves
the statement for fully decomposable $\tilde{\tilde p}$.
The general case follows from here.
Indeed, let $(B_i)\subset \z$, $i=1,\ldots,\deg \tilde{\tilde p}$, be the  coefficients
of $\tilde{\tilde p}$ (without loss of generality we can assume that the highest coefficient of $\tilde{\tilde p}$
is equal to unity).
The polynomial ring $k[B_i]$ is naturally embedded in the polynomial
ring $k[\la_i]$ where $\tilde{\tilde p}(u)=\prod_{i=1}^{\deg \tilde{\tilde p}}(u-\la_i)$. Since $\tilde{\tilde p}$ is fully decomposable
in $\z[\la_i]$, the statement
holds true for $\z[\la_i]$ and therefore for $\z[B_i]$, as the rings of scalars. This implies
the proposition via the evaluation homomorphism $\z[B_i]\to \z$.
\end{proof}
The meaning of Proposition \ref{pp'} is that the baxterization is essentially
independent on the particular polynomial annihilating the constant RE matrix $K$.
For instance, one may choose among them the  minimal, for  $p$.

Now let us look at the role of the parameter $\xi$.
Suppose that $\xi$ is invertible and define the polynomial
$p_\xi$ setting $p_\xi(u)=\xi^{-(n+1)}p(\xi u)$.
The correspondence $K\mapsto \xi K$ defines an algebra isomorphism $\K_p(x)\to \K_{p_\xi}(x)$, which
naturally extends to an isomorphism $\Ha_p(x)\to \Ha_{p_\xi}(x)$.
Therefore the image of a solution to (\ref{re}) is again a solution.
\begin{propn}
The isomorphism $\K_p(x)\to \K_{p_\xi}(x)$
sends $K_p(\xi,x)$ to $\xi^nK_{p_\xi}(1,x)$.
\end{propn}
\begin{proof}
Clear.
\end{proof}
\begin{remark}
Although the baxterization admits an arbitrary polynomial $p$, the latter is fixed in applications.
Thus one cannot get rid of the parameter $\xi$.
\end{remark}
\subsection{Unitarity and regularity of the principal solution}
In applications to integrable models, see e.g. \cite{S}, interesting are
those solutions to (\ref{re}) which pass through the unity at some $x_0\in k$;
they are called regular. Solutions are called unitary if
they obey the condition $K(x)K(x^{-1})=1$. Below we prove that the principal solution
(\ref{main}) is essentially unitary and regular.
\begin{propn}
\label{normalization}
The principal solution $K_p(\xi, x)$  satisfies the equality
\be
K_p(\xi, x)K_p(\xi, x^{-1})=\psi(\xi,x)\psi(\xi,x^{-1}),
\label{norm}
\ee
where
$\psi(\xi,x):=\frac{p(\xi x)}{\xi(x-x^{-1})}$.
\end{propn}
\begin{proof}
First of all note that $\psi$ is defined correctly. Indeed,
 $p(\xi x)$ is divisible by $\xi$, since $p(0)$
is divisible by the assumption.

Let us show that the element $K_p(\xi, x)K_p(\xi, x^{-1})$ actually belongs to $\z(x)$.
It is sufficient to check that assuming $\xi$ invertible, e.g. a nonzero
element from the field $k$ (cf. the proof of Theorem \ref{thm_main}).
Let $i$ be a non-negative integer, $i<n$.
The $K^{n-i}$-coefficient of the element $\xi(x-x^{-1})K_p(\xi, x)K_p(\xi, x^{-1})$ relative to the basis $\{K^j\}_{j=0}^n$ equals
\be
\xi(x-x^{-1})(b'\phi_i''+\phi_i'b''-\phi_n'\phi_i'' -\phi_n''\phi_i' +\sum_{\al=i}^{n} \phi_\al'\phi_{n+i-\al}''
+\sum_{j=0}^{n-1}A^{2n-j}_i\sum_{\al=0}^j \phi_\al'\phi_{j-\al}''),
\nn
\ee
where $b$ is defined by (\ref{b}).
We used the notation $a'=a(x)$ and $a''=a(x^{-1})$ for any $a\in \z(x)$.
Employing Lemmas \ref{dif} and \ref{sums'} (cf. the proof of Theorem \ref{thm_main}) we transform this expression to
\be
\xi(x-x^{-1})(b'\phi_i''+\phi_i'b''-\phi_n'\phi_i'' -\phi_n''\phi_i')+(\phi_{n+1}'\phi_{i}''-\phi_{i}'\phi_{n+1}'').
\nn
\ee
Now it is immediate to check that this expression vanishes.

To compute the value of $K_p(\xi, x)K_p(\xi, x^{-1})$, it is sufficient to consider the case when
$p(\mu)=0$ for some $\mu\in \z$ and $\xi$ is invertible in $\z$.
Under these assumptions let us apply Proposition \ref{pp'}
for  $\tilde p(u)=u-\mu$.
The homomorphism $\theta$ is identical on scalars
and therefore identical on $K_p(\xi,x)K_p(\xi,x^{-1})$.
We have
$$
K_p(\xi, x)K_p(\xi, x^{-1})=\theta\Bigl(K_p(\xi, x)K_p(\xi, x^{-1})\Bigr)= \Bigl(\psi(\xi,x)\frac{\xi x^{-1}-\mu}{\xi x-\mu}\Bigr)
\>\Bigl(\psi(\xi,x^{-1})\frac{\xi x-\mu}{\xi x^{-1}-\mu}\Bigr)
.
$$
This gives the right-hand side of (\ref{norm}).
\end{proof}
\begin{corollary}
Let $\z$ be a field. Then the element
$\frac{K_p(\xi,x)}{\psi(\xi,x)}$ is a unitary and regular
solution of the equation (\ref{re}).
\end{corollary}
\begin{proof}
Unitarity follows from Proposition (\ref{normalization}).
To check regularity, notice that both $K_p(\xi,x)$ and $\psi(\xi,x)$
have simple poles at $x=\pm1$ with equal residues $\frac{1}{2\xi}p(\pm \xi)$.
\end{proof}
\section{On uniqueness of the solutions found}
\label{s2T}
In this section we assume that $\z$ is an integral domain and denote by  $\Q\z$ its field of ratios.
We impose the condition $A_n=p(0)\not =0$ on the free term of the
polynomial $p$.
Without loss of generality, we will assume $A_n$ to be invertible in $\z$ (otherwise
we can pass to the localization $\z[A_n^{-1}]$).
Then $K$  is invertible in $\K_p\subset \Ha_p$, and for every integer $m$
we have
\be
\label{m-base}
K^{m+1}=A_0K^m+A_1K^{m-1}+\ldots+A_nK^{m-n}.
\ee
For each $m\in \Z$ the monomials $\{K^i\}_{i=m-n}^m$ form a $\z$-basis in $\K_p$.
The elements $\{e_{j;i}\}_{m-n\leq i<j\leq m}$ form a $\z$-basis in $\Lc_p$.
We introduce the dual basis $\{\al^{(m)}_{j;i}\}_{m-n\leq i<j\leq m}\subset \Lc_p^*$
to $\{e_{j;i}\}_{m-n\leq i<j\leq m}$ and rewrite
the system (\ref{system}) in the equivalent form
\be
(\al^{(m)}_{j;i}\circ Q_p)\bigl(K(x)\bigr)=0, \quad m-n\leq i<j \leq m
\label{system_m}.
\ee
This system is equivalent to the equation (\ref{re}).
\subsection{Necessary conditions for $K(x)$ to be a solution}
Equation (\ref{re}) has  been solved directly for $n=1$ and $n=2$.
In the present subsection we assume $n\geq 3$.

Let $K(x)$ be a solution to the equation (\ref{re}).
In terms of the basis $\{K^i\}_{i=m-n}^m\subset \K_p$, where $m$ is an integer,
the element $K(x)$ admits the decomposition
\be
K(x)&=&a_{m,0}K^{m}+a_{m,1}K^{m-1}+a_{m,2}K^{n-2}+\ldots +a_{m,n}K^{m-n}.
\label{K(x)}
\ee
Here $a_{m,k}$ are some coefficients from $\z(x)$.
In the present subsection  we establish
links among $a_{m_1,0}$ and $a_{m_2,0}$ for different $m_1$ and $m_2$. On this base, we
classify solutions to (\ref{re}) in Subsections \ref{subsecPS} and \ref{subsecSS}.

With any solution to (\ref{re}) we associate  the interval $I\subset [2,n]$ of integers defined as
$I=[2,n]$ if $a_{2,0}\not=0$ and $I=[3,n]$ otherwise.
\begin{lemma}
Let $K(x)\in \Ha_p$ be a solution to (\ref{re}).
Then there exists a set of constants
 $\{\xi_m\}_{m\in I}\subset \Q\z$ such that
\be
a_{m,1}&=&(\xi_mx-A_0)a_{m,0}
\label{1-0}
\ee
for all  $m\in I$.
\label{m,m-1}
\end{lemma}
\begin{proof}
For $m\geq 2$ the equation
$$
\bigl(\al^{(m)}_{m;m-1}\circ Q_p\bigr)\bigl(K(x)\bigr)=0,
$$
from the system (\ref{system_m}) reads
\be
(\frac{1}{y}-\frac{1}{x})A_0a_{m,0}'a_{m,0}''+
\frac{1}{y}a_{m,0}'a_{m,1}''
-\frac{1}{x}a_{m,1}'a_{m,0}''=0.
\label{aux_eq7}
\ee
Since $\z$ is an integral domain, so is the ring $\z[x]$. Passing to
its field of ratios, we find the general solution to (\ref{aux_eq7}).
That is either $a_{m,0}=0$ or
\be
a_{m,1}
=(\xi_m x -A_0)a_{m,0}
\nn
\ee
for some $\xi_m\in \Q\z$. Therefore the constant $\xi_m$ of interest does exist
for every $m\in[2,n]$ such that $a_{m,0}\not=0$.

Suppose $a_{m,0}=0$ for some $m\in [3,n]$ and consider the equation
$$
\bigl(\al^{(m)}_{m;m-2}\circ Q_p\bigr)\bigl(K(x)\bigr)=0
$$
from the system (\ref{system_m}).
Explicitly it reads
\be
(\frac{1}{y}-\frac{1}{x})\Bigl((A^2_0+A_{1})a_{m,0}'a_{m,0}''+
A_0(a_{m,1}''a_{m,0}'+a_{m,0}''a_{m,1}')
+a_{m,1}''a_{m,1}'\Bigr)
+\nn\\
+\frac{1}{y }a_{m,2}''a_{m,0}'
-\frac{1}{x}a_{m,2}'a_{m,0}''=0
\label{m,m-2}.
\ee

Since $\z$ has no zero devisors,
the equality $a_{m,0}=0$  for some $m\in [3,n]$ implies $a_{m,1}=0$. In this case we put $\xi_m=0$.
This completes the proof.
\end{proof}
\begin{lemma}
\label{lm_highest}
Let $K(x)\in \Ha_p(x)$ is a solution to (\ref{re}). Then
\be
\label{zeros}
\begin{array}{cclll}
a_{m-1,0}&=&a_{m,0}\xi_mx,& m\in I.
\end{array}
\ee
\end{lemma}
\begin{proof}
Let us re-expand $K(x)$ over the basis $\{K^i\}_{i=m-n-1}^{m-1}$:
\be
K(x)&=&(A_0a_{m,0}+a_{m,1})K^{m-1}+(A_{1}a_{m,0}+a_{m,2})K^{m-2}+\ldots+
(A_{n-1}a_{m,0}+a_{m,n})K^{m-n}+\nn\\
&+&
A_na_{m,0}K^{m-n-1}.
\ee
From this we get the following recurrent relations among $a_{m,i}$ and $a_{m-1,i}$:
\be
\left\{
\begin{array}{ccc}
a_{m,1}&=&a_{m-1,0}-A_0a_{m,0},\\
a_{m,2}&=&a_{m-1,1}-A_{1}a_{m,0},\\
a_{m,3}&=&a_{m-1,2}-A_{2}a_{m,0},\\
&\ldots&\\
a_{m,n}&=&a_{m-1,n-1}-A_{n-1}a_{m,0},\\
\end{array}
\right.
\label{rec}
\ee
which are valid for all integer $m$.
The first line of (\ref{rec})
together with Lemma \ref{m,m-1} gives
\be
\begin{array}{cccc}
a_{m,1}&=&a_{m-1,0}-A_0a_{m,0},\\
a_{m,1}&=&(\xi_mx-A_0)a_{m,0},&m\in I.\\
\end{array}
\nn
\ee
From this we obtain (\ref{zeros}).
\end{proof}
\begin{corollary}
If $a_{m,0}=0$ for some $m\in [3,n]$, then
 $a_{i,0}=0$ for all $i\in [2,m]$.
\label{cor_zeros}
\end{corollary}
\begin{lemma}
If $K(x)\in \Ha_p(x)$ is a solution to (\ref{re}), then
either $a_{m,0}=0$ for all $m\in [2,n-1]$ or $a_{m,0}\not =0$  for all  $m\in [1,n]$.
In the latter case, all the constants  $\xi_m$ from Lemma   \ref{m,m-1} are non-zero and equal.
\label{aux}
\end{lemma}
\begin{proof}
If $a_{n,0}=0$ then the lemma holds true because $a_{m,0}=0$ for all $m\in [2,n]$, by Corollary \ref{cor_zeros}.
So, assume that $a_{n,0}\not =0$.
Let the integer $\ell$ be the minimal from the interval $[1,n]$  such that $a_{m,0}\not=0$
for all $m\in [\ell,n]$.
To prove the lemma, it is sufficient to show  that either $\ell=n$ or $\ell=1$.
Suppose the opposite, that is $\ell\in [2,n-1]$.

Upon elementary transformations, the equation (\ref{m,m-2}) takes the form
\be
(\frac{1}{y }-\frac{1}{x})\Bigl(A_{1}+
\xi_m^2y x\Bigr)a_{m,0}'a_{m,0}''
+
\frac{1}{y }a_{m,2}''a_{m,0}'
-\frac{1}{x}a_{m,2}'a_{m,0}''=0.
\label{aux_eq2}
\ee
This holds for all $m\in[3,n]$. Since $\ell\in [2,n-1]$,  then
for  $m\in [\ell+1,n]$ the general solution to (\ref{aux_eq2}) is
\be
a_{m,2}
=(\xi_m^2x^2+\zt_m x-A_{1})a_{m,0},
\label{aux_eq3}
\ee
where $\zt_m\in \Q\z$ is an arbitrary constant.

Consider the equation
$$
\bigl(\al^{(m+1)}_{m;m-1}\circ Q_p\bigr)\bigl(K(x)\bigr)=0,
$$
assuming $m\in [\ell,n-1]$.
Explicitly, this equation reads
\be
(\frac{1}{y}-\frac{1}{x})\Bigl(A_0A_1 a_{m+1,0}'a_{m+1,0}''+A_1(a_{m+1,0}'a_{m+1,1}''+
a_{m+1,1}'a_{m+1,0}'')\Bigr)+\nn\\
\frac{1}{y}a_{m+1,1}'a_{m+1,2}''-\frac{1}{x}a_{m+1,2}'a_{m+1,1}''=0.
\label{aux_eq8}
\ee
Substitute in here the expressions for $a_{m+1,1}$ and $a_{m+1,2}$,
\be
a_{m+1,1}=(\xi_{m+1}x-A_0)a_{m+1,0}
,\quad a_{m+1,2}=(\xi_{m+1}^2x^2+\zt_{m+1} x-A_{1})a_{m+1,0}
\label{aux_eq9}
\ee
from (\ref{1-0}) and (\ref{aux_eq3}), respectively.
Then it is easy to check that (\ref{aux_eq8}) is fulfilled if and only if
$$\xi_{m+1}(\zt_{m+1} +\xi_{m+1}A_0)=0.$$
By assumption, $m\in [\ell,n-1]$ and therefore $a_{m,0}\not=0$. Then $\xi_{m+1}\not=0$, as follows from
(\ref{zeros}). Hence (\ref{aux_eq8}) is equivalent to
$\zt_{m+1}=-\xi_{m+1}A_0$.
Substitute this in the right equation of
(\ref{aux_eq9}) and get
$$a_{{m+1},2}=(\xi_{m+1}^2x^2-\xi_{m+1}A_0 x-A_{1})a_{m+1,0}.$$
Combine this with the equality $a_{m+1,2}=a_{m,1}-A_1a_{m+1,0}$ from (\ref{rec}) and get
$$
a_{m,1}=(\xi_{m+1}^2x^2-\xi_{m+1}A_0 x)a_{m+1,0}=(\xi_{m+1}x-A_0)\xi_{m+1}x a_{m+1,0}=(\xi_{m+1}x-A_0)a_{m,0}
$$
for all $m\in [\ell,n-1]$. In the rightmost equality we have used (\ref{zeros}).
Compare this equality with (\ref{1-0}). Since $a_{m,0}\not=0$ for $m\in [\ell,n-1]$,
this implies $\xi_{m}=\xi_{m+1}\not=0$. But then $a_{m-1,0}\not =0$, by the formula (\ref{zeros})
and $a_{\ell-1,0}\not =0$, in particular.
On the other hand, $a_{\ell-1,0} =0$, by definition of $\ell$.
This contradiction implies  either $\ell=n$ or $\ell=1$
(recall that we have supposed $\ell\in [2,n-1]$).
In the latter case, we have obtained $\xi_{m}=\xi_{n}\not=0$ for $m\in [2,n]$.
The proof is complete.
\end{proof}
\subsection{Principal solutions}
In the present subsection we assume $n\geq 3$.
Recall that a solution $K(x)$ is principal if $a_{m,0}\not=0$ in the expansion (\ref{K(x)}) for $m=n-1$ and hence for all $m\in[1,n]$,
by Lemma \ref{aux}. This lemma
also states that all the constants $\xi_m$ from Lemma \ref{m,m-1}
are equal to a non-zero element from $\Q\z$, which we denote
by $\xi$
\label{subsecPS}
\begin{lemma}
Let $K(x)\in \Ha_p$ be a principal solution to (\ref{re}).
Fix an integer $k\in [0,n-1]$. Then
\be
\begin{array}{cclll}
a_{m,k}
&=&\phi_{k}(\xi,x)a_{m,0},&m\in[k+1,n],
\end{array}
\label{k-s}
\ee
where $\phi_k$ are defined by (\ref{phi}).
\label{non-deg}
\end{lemma}
\begin{proof}
For $k=0$ the equality (\ref{k-s}) is true because $\phi_0=1$. So assume $k>0$.
Since $K(x)$ is principal, $a_{2,0}\not=0$ and $I=[2,n]$.
The first line of (\ref{rec})
together with Lemma \ref{m,m-1} gives
\be
\left\{
\begin{array}{lcccc}
a_{m,1}&=&a_{m-1,0}-A_0a_{m,0},\\
a_{m,1}&=&(\xi x-A_0)a_{m,0}&m=[2,n].
\end{array}
\right.
\nn
\ee
From this we get (\ref{k-s}) for $k=1$.

We will use induction on $k\in [1,n-1]$. The case $k=1$ is already checked.
Suppose the formula (\ref{k-s}) is proven for some $k\in [1,n-2]$.
For all $m=[2,n]$ the formula (\ref{rec}) gives
\be
\begin{array}{cclll}
a_{m,k+1}&=&a_{m-1,k}-A_{k}a_{m,0}
\end{array}
\nn
\ee
From this we get, using the  induction assumption,
\be
\begin{array}{cclll}
a_{m,k+1}&=&\phi_{k}a_{m-1,0}-A_{k}a_{m,0},
\end{array}
\nn
\ee
for
$m-1\in [k+1,n-1]$.
By (\ref{zeros}) we conclude that the equality $a_{m-1,0}=\xi x a_{m,0}$
is valid for all
$m\in[(k+1)+1,n]$.
Substituting this to the above equality we prove (\ref{k-s}),
according to the induction principle.
\end{proof}
\begin{corollary}
Fix $n\geq 2$. If $\z$ is a field and $p(0)\not=0$, then up to a factor
from $\z(x)$ the principal solution
to (\ref{re}) equals $K_p(\xi,x)$, where $\xi$ is an arbitrary  non-zero element
from $\z$
and $K_p(\xi,x)$ is given by (\ref{main}).
\label{unique}
\end{corollary}
\begin{proof}
Since $\z$ is a field,  we can put $a_{n,0}=1$
in the expansion (\ref{K(x)}).
Lemma \ref{non-deg} then uniquely determines the coefficients $a_{n,i}=\phi_i\in \z(x)$
for $i\in [2,n-1]$.
The free term $a_{n,n}$ is determined in the proof
of Theorem \ref{thm_main} and
has the form as in (\ref{main}).
This proves the statement.
\end{proof}

\subsection{Small solutions}
\label{subsecSS}
In this subsection we will clarify the structure of small solutions.
We know that for $n=1,2$ all non-trivial solutions are principal. Therefore we assume $n\geq 3$.
According to Lemma \ref{aux}, a small solution has $a_{m,0}=0$ for $m\in[2,n-1]$.

\begin{lemma}
\label{m-i}
Let $K(x)\in \Ha_p(x)$ be a small solution to (\ref{re}).
If $a_{m,0}=0$ for some $m\in [3,n]$,  then $a_{m,i}=0$ for all $i\in [0,m-2]$.
\end{lemma}
\begin{proof}
First of all, observe that for small solutions the interval $I$ coincides with  $[3,n]$,
since $a_{2,0}=0$.

We will use induction on $m$.
For $m=3$ the statement follows from Lemma \ref{m,m-1}.
Suppose that the statement  is proven for some  $m\in [3,n-1]$.
If $a_{m+1,0}=0$, we have  $a_{m+1,1}=0$, by Lemma  \ref{m,m-1}.
Then the first line in (\ref{rec}) implies $a_{m,0}=0$ (replace $m$ by $m+1$ there).
By the induction, we have  $a_{m,i}=0$ for $i\in [0,m-2]$.
Then (\ref{rec}) with  $m$ replaced by $m+1$ implies $a_{m+1,i}=0$ for $i\in [1,m-1]$.
Thus we have shown that $a_{m+1,0}=0$ implies $a_{m+1,i}=0$ for all $i\in [0,m-1]$.
By the induction principle, the lemma is proven.
\end{proof}
\begin{propn}
\label{pr_small}
Let $n\geq 3$ and $K(x)$ be a small solution to (\ref{re}). Then
$a_{n,i}=-A_{i-1} a_{n,0}$ for $i\in [1,n-2]$.
\end{propn}
\begin{proof}
Recall that a small solution has $a_{n-1,0}=0$ in the expansion (\ref{K(x)}).
If $n=3$ and $a_{2,0}=0$, then (\ref{rec}) implies $a_{3,1}=-A_0a_{3,0}$, as required.
Suppose $n\geq 4$. As follows from Lemma \ref{m-i}, $a_{n-1,i}=0$ for all $i\in [0,n-3]$.
Then the statement follows from (\ref{rec}), where $m$ should be set $m=n$.
\end{proof}
\begin{corollary}
\label{small}
Let $\z$ be an integral domain and $p(0)\not=0$. Then
a non-trivial small solution to the (\ref{re}) has the form
\be
K(x)=a_+K+a_0+a_-K^{-1},
\label{apm}
\ee
for some rational functions $a_\pm$ and $a_0$ with values in $\z$.
It may exist only for $n=3$.
\end{corollary}
\begin{proof}
It follows from Proposition \ref{pr_small} that
\be
K(x)&=&(K^n-A_0K^{n-1}-\ldots-A_{n-2}K-A_{n-1})a_{n,0}\nn\\
&+&
(a_{n,n-1}+A_{n-2}a_{n,0})K+(a_{n,n}+A_{n-1}a_{n,0})\nn\\
&=&a_{n,0}A_{n}K^{-1}+(a_{n,n-1}+A_{n-2}a_{n,0})K+(a_{n,n}+A_{n-1}a_{n,0}).\nn
\ee
Now put $a_+:=a_{n,n-1}+A_{n-2}a_{n,0}$, $a_-:=A_{n}a_{n,0}$, and $a_0:=a_{n,n}+A_{n-1}a_{n,0}$.
This proves the first part of the statement.

To complete the proof, let us pass to  the basis $\{K^i\}_{i=1-n}^1$.
It is straightforward to see that neither first nor second difference terms in
(\ref{re'}) contribute to $\Span\{e_{0;j}\}_{j=1-n}^{-3}$.
On the other hand, the element $K^2=\sum_{i=0}^nA_iK^{1-i}$ in the third term does.
From this we conclude that $a_+A_i=0$ for all $i\in [4,n]$. That is,
either $a_+=0$ or $A_i=0$ for all $i\in [4,n]$. But in the latter case $K$ is not invertible
if $n\geq 4$. Thus either $n=3$ or  $a_+=0$. But the latter option gives a trivial solution.
Indeed, the projection of the left-hand side of (\ref{re'}) to $\Span\{e_{0;-2}\}$
equals
$$
(x-\frac{1}{y}-y+\frac{1}{x})a_-'a_-''+(\frac{1}{y}-\frac{1}{x})(a_-'a_-''+A_3a_+''a_+'),
$$
modulo sign. Since $\z$ is an integral domain, $a_+=0$ implies $a_-=0$.
Therefore a non-trivial small solution (with invertible $K$) may exist only if
$n=3$. It is given by (\ref{small_sol}).
\end{proof}

\section{Applications}
\label{secA}
In the present section we consider some applications of the above results
to quantum groups. For a guide in the theory of quantum groups, the reader
is referred to \cite{D,FRT}.

Put $k=\C$.
Assume that $q$ is not a root of unity.
Let $(V,\rho)$ be the natural $n+1$-dimensional representation of
$U_q\bigl(gl(n+1)\bigr)$ and $\{e^i_j\}_{i,j=1}^{n+1}$ be the standard basis of matrix units.
The Jimbo R-matrix
\be
q\sum_{i=1}^{n+1} e^i_i\tp e^i_i+
\sum_{i,j=1\atop i\not=j}^{n+1} e^i_j\tp e^j_i+
\omega\sum_{i,j=1\atop i<j}^{n+1} e^i_i\tp e^j_j
\in \End(V^{\tp 2})
\label{gl}
\ee
satisfies the Hecke condition (\ref{Hecke}).
It equals to $P(\rho\tp \rho)(\Ru)$,
where $P$ is the permutation in $V\tp V$ and $\Ru$ is the universal R-matrix
of $U_q\bigl(gl(n+1)\bigr)$.

Denote by $\A$ the reflection equation
algebra associated with the natural $n+1$-dimensional representation of
$U_q\bigl(gl(n+1)\bigr)$, \cite{KS}. It is the quotient of the free algebra
$k\langle K^i_j\rangle$, $i,j=1,\ldots, n+1$ by the ideal of relations
(\ref{re0}).
The matrix $K$ is annihilated by the characteristic polynomial
$\sum_{i=0}^{n+1} (-1)^i\si_i K^{n+1-i}$, where
$\si_0=1$ and $\{\si_i\}_{i=1}^{n+1}$ belong to the center of $\A$,
\cite{PS}.
It is known that the center is isomorphic to the polynomial algebra
generated by $\si_i$.

Put $\z$ to be the center of $\A$.
The pair $R\in \End(V^{\tp 2})$, $K\in \End(V)\tp \A$
defines a homomorphism of the affine Hecke algebra $\Ha$
to $\End(V^{\tp 2})\tp \A$. Moreover, this homomorphism
factors through a homomorphism of the cyclotomic
Hecke algebra $\Ha_p$, where $p$ is the characteristic
polynomial for the RE matrix $K$.
The function $K(x)$ from Theorem \ref{thm_cubic} gives a spectral
dependant RE-matrix as a polynomial in $K$.

It is known that the matrix $(\rho\tp \id)(\Ru_{21}\Ru)\in \End(V)\tp
U_q\bigl(gl(n+1)\bigr)$ satisfies the reflection equation and
hence defines a homomorphism $\A\to U_q\bigl(gl(n+1)\bigr)$. This gives
 a homomorphism from $\Ha_p$ to $\End(V^{\tp 2})\tp
U_q\bigl(gl(n+1)\bigr)$.

Any representation of $\A$ in a vector space $W$
gives an RE matrix with coefficients in $\End(W)$.
An irreducible module $W$ defines a character of the center of $\A$.
Then the RE matrix $K$ satisfies the polynomial identity with numeric coefficients
obtained by evaluating the central elements $\si_i$ in $\End(W)$.
This polynomial may be distinct from the minimal polynomial for $K$.
According to Proposition \ref{pp'}, one can use for baxterization
any polynomial annihilating $K$; the results will differ
by a scalar factor.

Representations of $\A$
may be restricted from representations of $U_q\bigl(gl(n+1)\bigr)$, but not all.
For example, the algebra $\A$ has a bigger supply of one-dimensional representations.
According to \cite{M} (see also \cite{KSS}), a class of characters\footnote{There exists one more class of non-semisimple
RE
matrices, with quadratic minimal polynomial. The description of that class is more complicated, \cite{M}.} of $\A$
are representable by semisimple RE matrices with either two or three eigenvalues,
one being zero in the latter case:
\be
\begin{array}{c}
\begin{picture}(110,110)
\put(-2,0){\line(0,1){100}}
\put(0,0){\line(0,1){100}}
\put(100,0){\line(0,1){100}}
\put(102,0){\line(0,1){100}}
\put(0,20){\dashbox{2}(80,80){~}}
\put(25,45){\dashbox{2}(30,30){~}}
\put(6,90){$\scriptstyle{*}$}
\put(5,84){$\scriptstyle{\mu+\la}$}
\put(18,78){$\scriptstyle{*}$}
\put(31,65){$\scriptstyle{*}$}
\put(37,58){$\scriptstyle{\la}$}
\put(45,51){$\scriptstyle{*}$}
\put(71,90){$\scriptstyle{*}$}
\put(65,84){$\scriptstyle{u_i}$}
\put(59,78){$\scriptstyle{*}$}
\put(6,26){$\scriptstyle{*}$}
\put(11,33){$\scriptstyle{v_i}$}
\put(18,38){$\scriptstyle{*}$}
\end{picture}
\end{array}
.
\label{char}
\ee
Here the two dashed boxes may have arbitrary (even zero) size.
Non-zero entries are situated in the diagonals of the bigger
dashed box.
The symmetrically allocated numbers $u_i$ and $v_i$ in the skew diagonal of
the bigger box are subject to the condition $u_iv_i=-\la\mu$.
This RE matrix has eigenvalues $\mu$, $\la$, and $0$
(when the bigger dashed box is strictly less than the entire matrix).

The RE matrix (\ref{char}) satisfies the cubic polynomial equation
$K^3=A_0K^2+A_1K+A_2$
with $A_0=\la+\mu$, $A_1=-\la\mu$, and $A_2=0$.
Put $\z=\C$. Formula (\ref{cubic}) gives the following solutions:
\be
K(x)&=&K^2+(\xi x-\la-\mu)K-\frac{x\xi^2-(\la+\mu)\xi +\la\mu x^{-1}}{x-x^{-1}},
\quad \forall \xi\in \C,
\nn\\
K(x)&=&K^2-(\la+\mu)K+\frac{\zeta-\la\mu x^{-1}}{x-x^{-1}},
\quad \forall \zt\in \C,
\ee
corresponding to $\zt=0$ and $\xi=0$, respectively.

Assuming the matrix (\ref{char}) non-degenerate, we have
\be
K(x)&=&K+\frac{\zt+(\la+\mu)x^{-1}}{x-x^{-1}},
\quad \forall \zt\in \C,
\ee
as follows from (\ref{quadratic}).

\end{document}